\documentclass[12pt]{article}
\usepackage{amsmath,amsfonts,amsthm}

\def\hybrid{\topmargin 0pt      \oddsidemargin 0pt
        \headheight 0pt \headsep 0pt
        \voffset=-0.5cm
        \textwidth 6.5in        
        \textheight 9in         
        \marginparwidth 0.0in
        \parskip 5pt plus 1pt   \jot = 1.5ex}

\hybrid

\def\beq{\begin{equation}}
\def\eeq{\end{equation}}
\def\bea{\begin{eqnarray}}
\def\eea{\end{eqnarray}}
\def\p{\partial}
\def\G{\Gamma}
\def\g{\gamma}
\def\s{\sigma}

\def\a{\alpha}
\def\b{\beta}

\def\R{{\cal R}}
\def\A{{\cal A}}

\def\D{{\cal D}}

\def\J{{\cal J}}

\def\M{{\cal M}}

\def\SP{{\cal S}}

\def\res{{\rm res}}

\def\Log{{\rm Log}}
\def\re{{\rm Re\,}}
\def\im{{\rm Im\,}}

\def \matrix #1 {\left(\begin{array}{cc} #1 \end{array}\right)}

\newtheorem{theo}{Theorem}[section]
\newtheorem{prop}[theo]{Proposition}
\newtheorem{cor}[theo]{Corollary}
\newtheorem{lem}[theo]{Lemma}

\theoremstyle{definition}
\newtheorem{dfn}[theo]{Definition}
\newtheorem{rem}[theo]{Remark}

\def\bbZ{{\mathbb Z}}
\def\bbC{{\mathbb C}}

\def\bbR{\mathbb R}

\begin{document}

\title{Amoebas, Ronkin function and Monge-Amp\`ere measures of algebraic curves with
marked points.}
\author{I.Krichever
\thanks{Department of Mathematics, Columbia University; High School of Economics, Kharkevich Institute for Informations Transmission Problems, and Landau Institute for Theoretical Physics, Moscow, Russia.
Research is supported in part by Simons foundation and in part
by The Ministry of Education and Science of the Russian Federation (contract 02.740.11.5194).}}

\date{}

\maketitle

\begin{abstract} A generalization of the amoeba and the Ronkin function of a plane algebraic curve
for a pair of harmonic functions on an algebraic curve with punctures is proposed. Extremal properties of $M$-curves are proved and connected with the spectral theory of difference operators with positive coefficients.
\end{abstract}

\section{Introduction}

The amoeba $\A_f$ of a holomorphic function $f:(\bbC^*)^n\rightarrow \bbC$ (where $\bbC^*=\bbC\setminus{0}$) is, by definition, the image in $\bbR^n$ of the zero locus of $f$ under the mapping $\Log:(z_1,\ldots, z_n)\rightarrow (\log|z_1|,\ldots, \log|z_n|)$.
The terminology was introduced by Gelfand, Kapranov and Zelevinsky in \cite{gkz} and reflects
the geometric shape of typical amoebas, that is  a semianalytic closed subset of $\bbR^n$ with tentacle-like asymptotes going off to infinity. All connected components of the amoeba complement $\A_f^c=\bbR^n\setminus \A_f$ are convex. When $f$ is a Laurent polynomial, then, as shown in \cite{fpt}, there is a natural injective map from the set of connected components of $\A^c_f$ to the set of integer points of Newton
polytop $\Delta_f$ of $f$. This injective map is defined by the gradient $\nabla \R_f$ of, the so-called Ronkin function introduced in \cite{ronkin}:
\beq\label{ronkin_f}
\R_f(x)=\frac{1}{(2\pi i)^n}\int_{\Log^{-1}(x)}\frac{\log|\,f(z_1,\ldots,z_n)|
\,dz_1\cdots dz_n}{z_1\cdots z_n}
\eeq
The Ronkin function $\R_f(x)$ is convex. Recall, that each convex function $u$ defines the associated Monge-Amp\`ere measure $Mu$. If $u$ is a smooth convex function on  $\bbR^n$, then $Mu=\det({\rm Hess} (u))v$, where ${\rm Hess}(u)$ is the Hessian matrix and $v$ denotes the Lebesque measure on $\bbR^n$. If $u$ is convex but not necessary smooth $\nabla u$ can still be defined as a multifunction, and the Monge-Amp\`ere measure of $u$ for any Borel set $E$ is defined as $Mu(E)=v(\nabla u(E))$.

Since the Ronkin function is affine linear in a connected component of $\A_f^c$, the support of the associated Monge-Amp\`ere measure $\mu:=M\R_f$ is in $\A_f$. The latter and the lower bound for the Hessian of the Ronkin function for $n=2$ established in \cite{passare} give the upper bound for the area of two-dimensional amoebas in terms of the area of the Newton polygon:
\beq\label{amoebabound}
{\rm Area} (\A_f)\leq \pi^2 {\rm Area} (\Delta_f)
\eeq

Recently concepts of the amoebas, the Ronkin function and the associated Monge-Amp\`ere
measure for $n=2$ have attracted additional interest and have become a major tool in the studies
of topological types of real algebraic curves \cite{mikhalkin1,mikhalkin2} and in the study of random surfaces which arise as height functions of dimer configurations \cite{okounkov1,okounkov2,okounkov3}.

In \cite{mikhalkin2} it was proved that curves,
for which the upper bound (\ref{amoebabound}) is sharp are defined over $\bbR$ (up to a multiplication by a constant) and, furthermore, their real loci are isotopic to {\it simple Harnack curves} defined for the first time in \cite{mikhalkin1}. Simple Harnack curves form a particularly restrictive type of  Harnack curves shown to be topologically unique in [5] (and later unique up to algebraic deformations in \cite{okounkov2}. Recall, that classical Harnack curves are not unique even topologically -- in 1876 Harnack has proposed a construction of several series of maximal curves in the plane with various topological arrangement of ovals.

In the study of limiting shapes of random surfaces simple Harnack curves arise as spectral curves of
two-dimensional periodic difference operators with positive coefficients. The corresponding Ronkin function turned out to be the Legendre dual to a surface tension of the model. The extremal properties of Harnack curves were crucial for the construction in \cite{okounkov1,okounkov2,okounkov3} and lead to many probabilistic implications.

The main goal of the present paper is to extend the notions of the amoeba and the
Ronkin function of a plane algebraic curve to the case of an arbitrary smooth algebraic curve with marked points and to make the first steps in the study of their properties. More precisely, for each algebraic curve $\G$ with $n$ marked points $p_\a$ we define a family of amoebas parameterized by two {\it imaginary normalized} meromorphic differentials $d\zeta_1, d\zeta_2$ on $\G$ having simple poles at the marked points. The imaginary normalization means that all periods of the differentials are pure imaginary (the periods of the differential $id\zeta_a$ are real, i.e. $id\zeta_a$ is a real-normalized differential in the sense of definition introduced in \cite{kr-grush1}). Note, that the space of imaginary normalized differentials having simple poles at three points
has real dimension two. Therefore, the generalized amoeba corresponding to an algebraic curve with three punctures is unique up to a linear transformation of $\mathbb R^2$. This particular case is of special interest due to its connections with the spectral theory of integrable two-dimensional difference operators (see Section 5).

The imaginary (real) normalized differentials of the third kind per se are not new. They were probably known to Maxwell (the real part of such differential is a single valued harmonic function on $\G$ which is the potential of electromagnetic field on $\G$ created by charged particles at the marked points); they were used in the, so-called, light-cone string theory, and played a crucial role in joint works of S.~Novikov and the author on Laurent-Fourier theory on Riemann surfaces and on operator quantization of bosonic strings \cite{fourie}.

For the first time, real normalization as defining property of quasi-momentum differentials in the spectral theory of linear operators with quasi-periodic coefficients was introduced in \cite{kr-real} and \cite{kr-av} (where they were called absolute normalized). A notion of real-normalized meromorphic differentials is central in the Whitham perturbation theory of integrable systems. These real-normalized differentials  were used in  \cite{kr-grush1} and \cite{kr-grush2} for a new proof the Diaz' theorem on dimension of complete subvarieties of ${\cal M}_{g,n}$, and for the proof of the vanishing of a certain tautological class. In \cite{kr-arb} they were used for the proof of Arbarello's conjecture.

In the most general form the notions of the generalized amoeba and the Ronkin function can be defined for
{\it any pair of harmonic functions on an algebraic curve with punctures}. In the last section we present
this general setup and some examples associated with imaginary normalized differentials of the second kind.

{\bf Acknowledgments.} The author would like to thank Grisha Mikhalkin and Ovidiu Savin for the numerous valuable comments, clarifications, and suggestions.

\section{Amoebas of algebraic curves with two imaginary normalized differentials}
Let $(\G, p_1,\ldots,p_n)$ be a smooth genus $g$ algebraic curve with $n$ distinct marked points.
Non-degeneracy of the imaginary part of the Riemann matrix of $b$-periods of a basis of the normalized holomorphic differential on $\G$ implies that
\begin{lem} For a given set of $n$ real numbers $a=\{a_\a\in \bbR\}$ such that $\sum_{\a=1}^n a_\a=0$ there is a unique meromorphic differential $d\zeta_a$ on $\G$ which is holomorphic outside the marked points $p_\a$ where it has simple poles with residues $a_\a$,
$\res_{p_\a} d\zeta_a= a_\a$,
and such that all periods of $d\zeta_a$ are {\it pure imaginary}:
\beq\label{normalization}
\re \oint_{c} d\zeta_a=0,\ \ \forall c\in H_1 (\G,\bbZ).
\eeq
\end{lem}
Notice, that the imaginary normalization (\ref{normalization}) implies that the real part $x_a$ of the abelian integral $\zeta_a(p)=\int^p d\zeta_a+c$, i.e.
\beq\label{realpart}
x_a(p)=\re \left(\int^p d\zeta_a+c\right)
\eeq
is a {\it single-valued harmonic function} on the punctured Riemann surface
\beq\label{g0}
\G_0:=\G\setminus \{p_1,\ldots,p_n\}.
 \eeq
The zeros $q_s, \  s=1,\ldots, 2g+n-2$ (non necessary distinct) of $d\zeta_a$ are critical points of $x_a$. The level set $x_a^{-1}(f)$ for $f$ not in the set of critical values of $x_a$ ($f\notin \{x_a(q_s)\})$ is the union of smooth connected cycles on $\G$. If $f=f_s:=x_a(q_s)$ for some $s$, then $q_s$ is a self-intersection point of the level set $x_a^{-1}(f_s)$.

Let us fix a pair of $n$-tuples of real numbers $a_{j}=\{a_{\a,j}\},\, \sum_\a a_{\a,j}=0, j=1,2,$. Then for each smooth algebraic curve $\G$ with $n$ marked points $p_\a$ we have two associated imaginary normalized differential which for brevity will be denoted throughout the paper by $d\zeta_j:=d\zeta_{a_j}$.
\begin{dfn} The amoeba $\A_{\SP}\subset \bbR^2$ associated with the data $\SP=\{\G,p_\a,a_{\a,j}\}$
is the image of the map $\chi:\G_0\rightarrow \bbR^2, \, \chi(p)=(x_1(p),x_2(p))$, where $x_j(p)$ are harmonic functions on $\G_0$ defined by the imaginary normalized meromorphic differentials $d\zeta_j$.
\end{dfn}
\begin{rem}
The action of $GL_2(\bbR)$ on pairs of imaginary normalized differentials
\beq\label{ltrans}
d\zeta_j'=c_{1,j}\,d\zeta_1+c_{2,j}\,d\zeta_2 ,\ \ \{c_{ij}\}\in GL_2(\bbR)\,.
\eeq
corresponds to a linear transformation of the amoeba.

\begin{rem} If all periods of the differentials $d\zeta_j$ are integer multiple of $2\pi id_j^{-1}$, where $d_j$ are integers, then the functions $z_j=e^{d_j\zeta_j}$ are single-valued meromorphic functions on $\G$. In that case $\G$ is an irreducible component of the normalization of a (possibly singular) plane curve $(\widetilde \G, z_1,z_2)$. If $\widetilde \G$ is a stable singular curve of the irreducible type, then $\A_\SP$ coincides (up to rescaling of coordinates) with the amoeba of the plane curve $(\widetilde\G,z_1,z_2)$.

To some extend the proposed setup can be seen as a closer of a union of certain loci of degree $d$ plane (or more generally toric) curves as the degree $d\to\infty$. In order to clarify the latter statement it is instructive to consider the basic example of curves with three marked points $p_1,p_2,p_3$. Using if needed the linear  transformation (\ref{ltrans}) we may assume without loss of generality that the imaginary normalized differential $d\zeta_j$ has residue $1$ at the marked point $p_j$ and residue $-1$ at the point $p_3$. Consider the locus $\Sigma_{g,d}$ in $\M_{g,3}$ where the periods of the corresponding differentials are in $2\pi i d^{-1}\mathbb Z$. Then as it was mentioned above the function $z_j=e^{d\zeta_j}$
is a meromorphic function on $\G$ with pole of degree $d$ at the point $p_3$ and zero of degree $d$ at $p_j$.
Hence $\G$ is a normalization of a singular plane curve of degree $d$. The conditions on periods defining $\Sigma_{g,d}$
are independent on a choice of a basis of cycles and in each fixed basis are rationality type condition. Therefore, the union of $\Sigma_{g,d}$ for all $d$ is dense in $\M_{g,3}$.
\end{rem}

\end{rem}
Our first goal is to describe explicitly  the locus of the critical points of the amoeba map:
\begin{lem} \label{critical} The locus $\g\subset \G_0$ of the critical points of the amoeba map $\xi:\G_0\rightarrow \bbR^2$ is a union of the locus $\g_0$, where the function
\beq\label{R}
R(p)=\frac{d\zeta_1}{d\zeta_2}
\eeq
is real, i.e. $\g_0:=\{p\in \g_0| \, \im R(p)=0\}$, and the finite set (possibly empty) of the common zeros of the differentials $d\zeta_j$.
\end{lem}
\begin{rem}
A pair of meromorphic differentials on a smooth algebraic curve defines a map of the curve to the two-dimensional complex projective space, $p\in \G\rightarrow (d\zeta_1(p):d\zeta_2(p))\in \mathbb{CP}^2$. For plane curves $f(z_1,z_2)=0$ where $d\zeta_j=d\ln z_j$ this map is known as the logarithmic Gauss map. The characterization of the critical points of the amoeba map for plane curves as the preimage under the logarithmic Gauss map of $\mathbb{RP}^2\subset \mathbb{CP}^2$ was obtained in \cite{mikhalkin1}.
\end{rem}
\noindent{\it Proof.}
The statement of the Lemma is a direct corollary of the following simple computations motivated by the computations of the Hessian of the Ronkin function in \cite{passare}.

Let $p$ be a regular point of the map $\chi$, i.e. at $p$ the form $dx_1\wedge dx_2$ is non-degenerate. Then in the neighborhood of $p$ the functions $(x_1(p),x_2(p))$ define a system of real coordinates and we can write
\beq\label{calc1}
\zeta_1=x_1+iy_1(x_1,x_2),\ \ \zeta_2=x_2+iy_2(x_1,x_2)
\eeq
The imaginary parts $y_j$ of the abelian integrals $\zeta_j$ are multivalued functions
with periods of $d\zeta_j$.

Taking the full external differential of equations (\ref{calc1}) we get
\beq\label{calc2}
d\zeta_1=(1+i\p_1y_1)\,dx_1+(i\p_2y_1)\,dx_2,\ d\zeta_2=(i\p_1y_2)\,dx_1+(1+i\p_2y_2)\, dx_2
\eeq
The ratio $R=d\zeta_1/d\zeta_2$ of two meromorphic differentials is a meromorphic function on $\G$. Hence,
\beq\label{calc3}
(i\p_1y_2)R=1+i\p_1y_1, \ \ (1+i\p_2y_2)R=i\p_2y_1
\eeq
From the first equation in (\ref{calc3}) we get
\beq\label{calc4}
\p_1y_2=-\frac{1}{\im \,R}\ ,\ \ \p_1y_1=-\frac{\re\, R}{\im \,R}
\eeq
The second equation in (\ref{calc3}) implies
\beq\label{calc5}
\p_2y_2=\frac{\re\, R}{\im \,R} ,\ \ \p_2 y_1=\frac{|\,R|^2}{\im\, R}
\eeq
From (\ref{calc2}) and equations (\ref{calc4},\ref{calc5}) it easy follows:
\beq\label{dx12}
d\zeta_1\wedge d\bar\zeta_2=(2+2i\p_1y_1)\,dx_1\wedge dx_2
\eeq
Hence,
\beq\label{vol}
4dx_1\wedge dx_2=d\zeta_1\wedge d\bar\zeta_2+d\bar \zeta_1\wedge d\zeta_2=-2i (\im R)\, d\zeta_1\wedge d\bar\zeta_1=-2i\left(\frac {\im R}{|\,R|^2}\right)\,d\zeta_2\wedge d\bar\zeta_2\,.
\eeq
For any meromorphic form $\omega$ the two-form $i d\omega\wedge d\bar \omega$ is finite and positive  except at poles and zeros of the differential.
Therefore, equation (\ref{vol}) implies that {\it non-isolated} critical points of $\chi$ is the locus $\g_0$
where $R$ is real. From (\ref{vol}) it follows also that isolated critical points of $\chi$ are
{\it common zeros} of $d\zeta_1$ and $d\zeta_2$.

\medskip

Note for further use, that equation (\ref{vol}) allows to identify the closed subsets $\G^+$ and $\G^-$  of $\G_0$, where the jacobian $J(\chi)$ of the map $\chi$ is non-negative and nonpositive, respectively, with the loci
\beq\label{orient}
\G^\pm:=\{p\in\G^\pm| \mp\im R(p)\geq0\}.
\eeq
The intersection $\g_0=\G^+\cap \G^-$ of these sets is the locus of non-isolated critical points of the amoeba map. It is a graph on $\G$ whose vertices are zeros of the differential $dR$ and non simple poles of $R$ (corresponding to multiple zeros of $d\zeta_2$). Generically, edges of the graph are preimages of the folding lines of the amoeba map.

The following statement shows that the geometric shape of a generalized amoeba is the same as that for
the typical amoebas of plane curves:
\begin{prop}\label{convexity}
All connected components of the complement $\A_\SP^c$ are convex. There are $n$ unbounded
components separated by tentacle-like asymptotes of the amoeba.
\end{prop}
\noindent{\it Proof.}
Convexity (local) of the amoeba compliment in the neighborhood of any point of its boundary
is an easy corollary of the maximum principle for harmonic functions $x_i$ or their linear combinations. The proof of the global convexity of each connected component requires additional arguments. In the theory of conventional amoebas there are few different types of such arguments. One of the simplest is the use of convexity of the Ronkin function (see \cite{passare} or \cite{mikhalkin1}). Identically the same arguments using the generalized Ronkin function (which will be introduced in the next section) prove the first statement of the Proposition.

The second statement of the Proposition is obvious when the asymptotic directions
\beq\label{phi}
\varphi_\a:={\rm Arg}\,(a_{\a,1}+ia_{\a,\,2})
\eeq
of tentacles of the amoeba are all distinct. Notice, that for a given curve $\G$ with fixed marked points $p_\a$ the differentials $d\zeta_j$ depend continuously on the parameters $a_{\a,\,j}$. Hence, the intersection of the amoeba with any compact set in $\bbR^2$ also depends continuously on $a_\a$. That implies that the second statement of the Proposition holds in the general case, as well.

\begin{rem} It would be interesting to find an upper bound $\nu(g,n)$ for the number of connected
components of the generalized amoebas corresponding to smooth genus $g$ algebraic curves with $n$ marked points.
\end{rem}

\section{The generalized Ronkin function}

For each $x=(x_1,x_2)\in \bbR^2$ introduce the closed subsets $\G_x^{\pm}\subset \G_x\subset \G$:
\beq\label{gx}
\G_x:=\left\{p\in \G_x| \,x_1(p)\leq x_1,\, x_2(p)\leq x_2\right\}, \ \ \G_x^{\pm}=\G_x\cap \G^{\pm}.
\eeq
\begin{dfn} The generalized Ronkin function $\rho=\rho_\SP$ associated with a smooth algebraic curve with two imaginary normalized differentials ($\SP=\{\G,p_\a,a_{\a,j}\}$) is given by the formula:
\beq\label{ronkin}
\rho\,(x)=\frac {1}{8\pi i}\int\int_{\G_x}{\rm sgn}(\im R(p))\,(d\zeta_1\wedge d\bar\zeta_2-
d\bar\zeta_1\wedge d\zeta_2).
\eeq
\end{dfn}
\begin{theo}\label{Hess} The generalized Ronkin function $\rho\,(x)$ given by (\ref{ronkin}) is a convex function on $\bbR^2$. It is affine linear in each connected component of $\A_{\SP}^c$. It is smooth at the regular points of the amoeba, i.e. outside of the set $F$ of critical values of $\chi$, and furthermore at $x\in \A_\SP\setminus F$
\beq\label{hessian}
{\rm Hess}\, \rho\,(x)=\frac1{2\pi}\sum_{p\,\in \,\chi^{-1}(x)}\frac{1}{|\,\im R(p)|}\left(\begin{array}{cc} 1 &\re R(p)\\
\re R(p)& |\,R(p)|^{\,2}
\end{array} \right).
\eeq
\end{theo}
\begin{rem} Formula (\ref{hessian}) in the case of plane curves coincides with the formula for the hessian of Ronkin function $\R_f$ obtained in $\cite{passare}$. Hence, in the case of plane curve the function $\rho$ coincides with the conventional Ronkin function up to an affine linear form.
\end{rem}
\noindent{\it Proof.} Let $x^0$ be in a complement of the amoeba. Then the boundary of the domain $\G_{x^0}$ is a disjoint
union of the cycles $l_j(x^0):=x_j^{-1}(x_j^{0})\,\cap \,\G_{x^0}, \, j=1,2$. For a generic $x^0\in \A_\SP^c$ these cycles are smooth. It will be always assumed that segments $l_j^\pm (x^0)=l_j(x^0)\cap \G^\pm$ are oriented such that $dx_k,\ k\neq j,$ is positive with respect to the orientation. Then, $l_j(x^0)=l_j^+(x^0)-l_j^-(x^0)$.

Equation (\ref{dx12}) implies that the generalized Ronkin function can be represented in the form:
\beq\label{ronkin2}
2\pi \rho\,(x)=\int\int_{\G_x^+} (\p_1y_1)dx_1\wedge dx_2-\int\int_{\G_x^-} (\p_1y_1)dx_1\wedge dx_2
\eeq
Hence,
\beq\label{grad2}
2\pi\p_2 \rho\,(x^0)=\int_{l_1^+(x^0)} (\p_1y_1)dx_1-\int_{l_1^-(x^0)} (\p_1y_1)dx_1=-i\oint_{l_1(x^0)}d\zeta_1
\eeq
From (\ref{calc4}) and (\ref{calc5}) it follows that $\p_1 y_1=-\p_2y_2$. Therefore,
\beq\label{grad1}
2\pi \p_1 \rho_{\SP}(x^0)=i\oint_{l_2(x^0)}d\zeta_2
\eeq
The periods of the differentials in the right hand sides of (\ref{grad2},\ref{grad1}) are constant in each of the connected component of $\A_\SP^c$. Therefore, $\rho$ is affine linear in each of these components.

Consider now a regular point $x$ of the amoeba, i.e. $x\in \A_\SP\setminus F$, where $F$ is the set of critical values of $\chi$, then from (\ref{ronkin2}) we get the equation:
\beq\label{grad3}
2\pi \p_1 \rho\,(x)=-\sum_{p\in \,\chi^{-1}(x)\cap \G^+} y_2(p)+\sum_{p\in \,\chi^{-1}(x)\cap \G^-} y_2(p)+\Pi_2(x),
\eeq
where $\Pi_2$ is the imaginary part of the period of $d\zeta_2$ (over connected components of the level set $l_2(x)$ which do not intersect $l_1(x)$).
Similarly, we get
\beq\label{grad4}
2\pi \p_2 \rho\,(x)=\sum_{p\in \,\chi^{-1}(x)\cap \G^+} y_1(p)-\sum_{p\in \,\chi^{-1}(x)\cap \G^-} y_1(p)+\Pi_1(x),
\eeq
Notice that the periods $\Pi_1$ and $\Pi_2$ are locally constant. Therefore, taking the derivatives of (\ref{grad3}) and (\ref{grad4}) and using the equations (\ref{calc4}) and (\ref{calc5}) we get equation (\ref{hessian}) for the Hessian of $\rho$ at the regular points of the amoeba.

Each term in the sum at the righthand side in (\ref{hessian}) is a symmetric matrix with determinant 1
and with positive diagonal elements. Hence, each term of the sum is positive definite. Therefore, ${\rm Hess}\, \rho$ is positive definite in the interior of the amoeba.

Let $\sigma$ be a subset of the points of the critical set $F$ whose preimages are not finite. If the primage of $x\in \s$ is not finite, then it is a closed cycle of the graph $\g_0$. Generically $\s$ is empty, but in any case it is a finite set.  From equations (\ref{grad3},\ref{grad4}) it easy follows that $\nabla \rho$ extends continuously  at $x\in F\setminus \s$. That implies the convexity of the generalized Ronkin function $\rho$. The theorem is proved.

Notice, that the associated Monge-Amp\'ere measure has nonzero point measure at each $x\in \s$. Indeed, on the cycle $c=\chi^{-1}(x)\in \g$ we have
\beq\label{yy}
\frac{dy_2}{dy_1}(t)=\frac{d\zeta_2}{d\zeta_1}(t), \ t\in c
\eeq
Hence, the latter function is not a constant. Therefore, the convex hull of the line $(-y_2(t),y_1(t))\in \bbR^2$ has non-trivial measure.

\medskip
From the proof of Theorem \ref{Hess} it follows that the image of the gradient map
$\nabla \rho\,: \bbR^2\rightarrow \bbR^2$ is a polygon $\Delta_\SP$
that is the convex hull of the set of points $v_\a=(v_{\a,1},v_{\a,2})\in \bbR^2$, which are the images under $\nabla \rho\,$ of the unbounded components of $\A_{\SP}$.

We describe first the polygon $\Delta_\SP$ in the case when asymptotes $\varphi_a$ given by (\ref{phi}) of the tentacles the amoeba are all distinct. Using if needed a linear transformation (\ref{ltrans}) we may assume also that $\a_{\a,j}\neq 0$.
The periods (\ref{grad2},\ref{grad1}) defining the image under $\nabla \rho\,$ of unbounded components are sums of periods over small cycles around subsets of marked points. The latter are equal (up to a sign) to the corresponding residues of the differential. It is straightforward to verify that the coordinates of the vertices of $\Delta_\SP$ equal
\beq\label{pol1}
{\rm If} \pm a_{\a,2}>0, \ \ {\rm then} \ \
v_{\a,1}=\pm\sum_{\b\in I_\a^\pm} a_{\b,2}, \ {\rm where}\ \ I_\a^{\pm}:=\{\b\in I_\a^\pm| \ \pm a_{\b,2}>0 , \pm (\varphi_\b-\varphi_\a)>0\}\,.
\eeq
\beq\label{pol2}
{\rm If} \pm a_{\a,1}>0, \ \ {\rm then} \ \
v_{\a,2}=\mp\sum_{\b\in J_\a^\pm} a_{\b,1}, \ {\rm where}\ \ J_\a^{\pm}:=\{\b\in J_\a^\pm| \ \mp a_{\b,2}>0 , \mp (\varphi_\b-\varphi_\a) >0\}\,.
\eeq
It is easy to see that the formulae (\ref{pol1}) and (\ref{pol2}) extend continuously to the general case when some of $\varphi_{\a}$ might coincide.
\begin{cor} The inequality
\beq\label{areabound}
{\rm Area} (\A_\SP)\leq \pi^2 {\rm Area} (\Delta_\SP)
\eeq
holds.
\end{cor}
\noindent
The proof of (\ref{areabound}) is identical to that in \cite{passare}.
Namely, as shown in \cite{passare}, for $(2\times 2)$  symmetric positive definite matrices $M_1$ and $M_2$ the inequality
\beq\label{passare1}
\sqrt{\det\, (M_1+M_2)}\geq \sqrt{\det M_1}+\sqrt{\det M_2}
\eeq
holds, with equality precisely if $M_1$ and $M_2$ are real multiples of each other.
The sum in the righthand side of (\ref{hessian}) contains at least two terms. Hence, $\det {\rm Hess}_\SP (x) \geq \pi^{-2}$ at the regular points of amoeba, $x\in \A_\SP\setminus F$. As shown above, the image under the gradient map of the regular points of the amoeba is $\Delta_\SP$ minus some of its vertexes. Combining the last two facts one gets
(\ref{areabound}).

\section{Amoebas of $M$-curves}

At the beginning of this section we show that the bound (\ref{areabound}) is sharp, and is achieved for certain pairs of imaginary normalized differentials on the, so-called, $M$-curves.

Recall that a smooth genus $g$ algebraic curve $\G$ with antiholomorphic involution $\tau:\G\rightarrow \G$ is $M$-curve if $\tau$ has maximal possible number of fixed cycles. By Harnack inequality, that number is $(g+1)$. The anti-involution $\tau$ on an $M$-curve is always of the separating type, i.e.
$\G$ is a union of two closed subsets $\G^\pm$ such that $\tau:\G^+\rightarrow \G^-$ and
$\G^+\cap \G^-=\cup_{j=0}^g A_j$, where $A_j$ are fixed ovals of $\tau$.

\begin{dfn} The set of data $(\G,p_\a, a_{\a,j})$ is called Harnack: (i) $\G$ is a $M$-curve; (ii)
the marked points $p_{\a}$ are on one of the fixed ovals $A_0,\ldots A_g$ of the antiinvolution
$\tau$, say $p_\a\in A_0$;
(iii) The cyclic order $p_\a$ along the cycle $A_0$ coincides with counterclockwise order of the vertices
of the polygon $\Delta_\SP$.
\end{dfn}
The following statement describes explicitly the zero divisor of an imaginary normalized differential on a $M$-curve under certain assumptions on positions of its poles and signs of its residues:
\begin{lem}\label{zeros}
Let $d\zeta$ be an imaginary normalized differential on a $M$-curve $\G$ with poles $p_\a$ on one of the fixed ovals of the involution, $p_a\in A_0$. If signs of the residues of $d\zeta$ at $p_\a$ are changed only twice with respect to cyclic order along the oval, then $d\zeta$
has no zeros outside of the fixed ovals of the anti-involution $\tau$. Furthermore: (i) all zeros of $d\zeta$ are simple, (ii) there is exactly $2$ zeros of $d\zeta$ on each of $g$ fixed ovals of $\tau$ that do not contain poles of $d\zeta$; (iii) there is exactly one zero of $d\zeta$ on each segment of $A_0$ between consecutive poles with the same sign of the residues.
\end{lem}
\noindent{\it Proof.} Notice, that the imaginary normalization condition is invariant with respect to
an antiholomorphic involution, i.e. $d\zeta=\tau^*(d\bar\zeta)$. Hence, $d\zeta$ is real on the fixed ovals of $\tau$. Then, the imaginary normalization implies that the period of $d\zeta$ along  $A_j$ which do not contain the marked points vanishes. Hence, on each of these ovals there are at least
2 distinct zeros of $d\zeta$.

Let $p_\a$ and $p_{\a+1}$ be two consecutive (along $A_0$) marked points with the same sign of the residues. Then,
$d\zeta$ has at least one zero on the segment of $A_0$ between these points.
By the second assumption of the lemma the number of consecutive pairs of the marked points with the same signs of the residues is at least $n-2$. The zero divisor of $d\zeta$ is of degree $2g+n-2$. Hence, all the previous bounds on the number of distinct zeros are sharp and we have accounted all the zeros of $d\zeta$. The lemma is proved.

\begin{lem} \label{diff}
The amoeba map defined by the imaginary normalized differentials $d\zeta_j$ associated with
Harnack data restricted to $\G^+\subset \G$ is a diffeomorphism of $\G_0^+:=\G^+\setminus
\{p_\a\}$ with $\A_{\SP}$.
\end{lem}
\noindent{\it Proof.} Let us show first, that for the Harnack pairs the set $\g$ of critical points of
the amoeba map coincides with the locus $A_0,\ldots, A_g$ of fixed points of the anti-involution $\tau$. Indeed, suppose not, then as follows from Lemma (\ref{critical}) there exists a point $p\notin A_j$ such that $R(p)=r$ is real. By definition of $R$ the differential $d\zeta=d\zeta_1-rd\zeta_2$
has zero at the point $p$. Notice now, that if $d\zeta_1, d\zeta_2$ is a Harnack pair of the imaginary normalized differentials, then the differential  $d\zeta$ satisfies the assumptions of Lemma \ref{diff}, and can not vanish outside the fixed ovals of the involution.

If the locus of $F$ of critical points of the amoeba map coincides with the boundary of $\A_\SP$ then the function $x_2$ restricted to a segments of the level line $x_1^{-1}(a)$ in $\G_0^+$ is monotonic. The lemma is proved.

The same arguments show that
\begin{lem} For Harnack data the gradient $\rho$ restricted to $\A\setminus F$ is one-to-one map onto the interior of the polygon $\Delta_S$ with $g$ points removed. The latter are images under $\rho$ of the fixed ovals $A_1,\ldots,A_g$.
\end{lem}

Our next goal is to show that as in the case of plane curves the extremal property of the area characterizes the Harnack pairs.

\begin{theo} The equality
\beq\label{areaeq}{\rm Area} (\A_\SP)=\pi^2 {\rm Area} (\Delta_\SP)\,.
\eeq
holds if and only if
$\A_\SP$ is associated with a Harnack pair of the imaginary normalized differentials.
\end{theo}
\noindent
{\it Proof.} The key arguments in the proof are almost identical to that in \cite{mikhalkin2}. First notice, that if the equality (\ref{areaeq}) holds then the sum in (\ref{hessian}) contains two equal terms. Hence on the preimage
$\chi^{-1}(\A\setminus F$ there is an involution $\tau$ and
\beq\label{involution}
R(z)=\bar R(\tau(z)).
\eeq
From the latter it follows that $\tau$ is an anti-holomorphic involution of $\chi^{-1}(\A\setminus F)$. Another corollary of (\ref{areaeq}) is that there is no points of $F$ having non-zero point mass, i.e. all points of $\A$ has finite number of preimages. From that it easy follows that each point of $F$ has only one preimage. Hence, $\tau$ extends to an anti-holomorphic involution of $\G$ and preimages of bounded connected components of $F$ are fixed ovals of the involution. Preimages of unbounded connected components of F are segments of $A_0$ between consecutive poles $p_\a$

\section{Two-dimensional difference operators.}

For completeness, in this section we briefly outline some relevant facts from the spectral theory of difference operators. In the framework of the spectral theory of two-dimensional difference operators
\beq\label{L}
(L\psi)_{n,m}=\psi_{n+1,\,m}+\psi_{n,\,m+1}+u_{n,\,m} \psi_{n,\,m}
\eeq
acting on functions of the discrete variables $(n,m)\in \bbZ^2$ (see \cite{kr-diff,okounkov2}) the plane curves of degree $d$ arise as the spectral curves of $d$-periodic operators, $v_{n,\,m}=v_{n+d,\,m}=v_{n,\,m+d}$. The points of these curve parameterize, the so-called, Bloch solution of the equation $L\psi=0$, i.e. the solutions that are eigenfunctions of the monodromy operators
\beq\label{monod}
\psi_{n+d,\,m}=z_1\psi_{n,\,m},\ \ \psi_{n,\,m+d}=z_2 \psi_{n,\,m}.
\eeq
In \cite{kr-diff} an explicit solution of the {\it inverse spectral transform} in terms of Riemann
theta-function was obtained for {\it integrable} difference operators. In the periodic case the
integrable difference operators (related to solutions of the bilinear discrete Hirota equations
(BDHE)) are singled out by the condition that the affine spectral curve is compactified
by {\it three infinite points}.

Recall, briefly the solution of the inverse problem. Consider a smooth genus $g$ algebraic curve
$\G$  with three marked points $p_1,p_2,p_3$. Then, a choice of a basis of cycles $\{A_i,B_i\}$ on $\G$ with the canonical matrix of intersections: $A_i\cap B_j=\delta_{i,\,j},\ A_i\cap A_j=B_i\cap B_j=0$ allows one to introduce the basis of normalized holomorphic differentials $\omega_i, \oint_{A_j}\omega_i=\delta_{i,j}$ and the matrix $B$ of its periods, $B_{i,j}=\oint_{B_j}\omega_i$.
The latter is a symmetric matrix with positive definite imaginary part. It defines the Jacobian
of the curve $\J(\G):=\mathbb C^g/\Lambda$, where  the lattice $\Lambda$ is generated by the basis vectors $e_m\in \mathbb C^g$ and the column-vectors of $B$.
The Riemann theta-function $\theta(z)=\theta(z|B)$ corresponding to $B$
is given by the formula
\beq\label{teta1}
\theta(z)=\sum_{m\in \mathbb{Z}^g} e^{2\pi i(z,m)+\pi i(Bm,m)},\ \
(z,m)=m_1z_1+\cdots+m_gz_g\, .
\eeq
The multi-valued correspondence $A(p)\in \bbC^g, A_k(p)=\int^p\omega_k$ descents to the so-called Abel map $A:\G\rightarrow \J(\G)$.

Introduce also the  normalized abelian integral $d\Omega_1 (d\Omega_2)$ of the third kind having simple pole with residue $1$ at $p_1(p_2)$ and simple pole with residue $-1$ at $p_3$. Their abelian integrals $\Omega_j(p)$ are multivalued functions on $\G$.

\begin{theo}(\cite{kr-diff}) Let $p$ be a point of $\G$. Then for any $g$-dimensional vector
$Z$ the function
\beq\label{pd}
\psi(m,n)=\frac{\theta(A(p)+mU+nV+Z)}{\theta(mU+nV+Z)}\, e^{m\Omega_1(p)+n\Omega_2(p)},
\eeq
where $U=A_(p_1)-A(p_3),\ V=A(p_2)-A(p_3)$, satisfies the difference equation
\beq\label{laxdd}
\psi(m,n+1)=\psi(m+1,n)+u(m,n)\psi(m,n)
\eeq
with
\beq\label{ud}
u(m,n)=\frac{\tau(m+1,n+1)\tau(m,n)}{\tau(m+1,n)\tau(m,n+1)}
\eeq
where
\beq\label{tau}
\tau(m,n)=c_1^mc_2^n \,\theta(mU+nV+Z)
\eeq

\end{theo}
A few remarks: (a) In general the coefficients of the difference equation are complex qusiperiodic function (possibly singular) of the variables $(m,n)$. (b) If the vectors $U$ and $V$ are $d$-periodic points of the Jacobain, then $u(m,n)$ are $d$-periodic. (The latter condition is equivalent to the condition that on $\G$ there exist function $z_1$ and $z_2$ having pole of order $d$ at $p_3$ and zeros of order $d$ at $p_1$ and $p_2$, respectively). (c) If on $\G$ there is an antiholomorphic involution and if the marked points $p_j$ are fixed by the involution, then $u$ is real for real $Z$ (but still might be singular). Finally, (d)
\begin{lem}\label{nonsing}
If $\G$ is $M$-curve, and the marked points are on one of the fixed ovals of the involution, then for all real $Z$ the coefficients of the difference equation are non-singular for all $(n,m)$.
\end{lem}
The proof of the lemma is standard in the finite-gap theory. Namely, the corresponding discrete Baker-Akhiezer function
$\psi(m,n,p)$ is a unique meromorphic  function on $\G$ having pole of order $n+m$ at $p_3$, zero of order $m$ at $p_1$ and zero of order $n$ at $p_2$, and having in addition at most simple poles at points $\g_1,\ldots,\g_g$ (if the latter are distinct), i.e. $\psi(m,n,p)\in H^0((m+n)p_3-mp_1-np_2+D),\ D=\g_1+\cdots+\g_g$. In addition to zeros at $p_1,p_2$,
this function has $g$ zeros $\g_s(n,m)$. The coefficient $u(m,n)$ is singular iff one of these zeros coincides with $p_3$.

Let $\G$ be an $M$-curve, and the marked points be on one of the fixed ovals of the involution. If points $\g_s$ are chosen such that $\g_s\in A_s$, then from the uniqueness of the Baker-Akhiezer function it follows that $\bar\psi=\tau^*\psi$, i.e. $\psi$ is real on each of the fixed ovals of the antiinvolution $\tau$. On each of this ovals it has pole, hence it must have a zero. The total number of these zeros $\g_s(m,n)$ is $g$. Hence, none of them can coincide with $p_1\in A_0$. The lemma is proved.

Under the gauge transformation $\Psi(m,n,p)=(-1)^n\tau(m,n)\psi(n,m)$, equation (\ref{laxdd}) takes the following form
\beq\label{laxdd1}
\tau(m+1,n)\Psi(m,n+1)+\tau(n,n+1)\Psi(m+1,n)+\tau(m+1,n+1)\Psi(m,n)=0
\eeq
\begin{cor}Under the assumption of Lemma \ref{nonsing} the coefficients of equation \ref{laxdd1} are real positive numbers.
\end{cor}
In order to prove the corollary it is enough to consider a degeneration of
$\G$ to the rational curve with three fixed points. Such a degeneration can be performed within the space of $M$-curves,
therefore along the deformation the coefficients remain non-zero, and the operator corresponding to a rational curve has constant positive coefficients.

\section{Concluding remarks: Amoebas versus jellyfish and other creatures.}

In the most general form the amoeba map $\chi: \G_0\rightarrow \bbR^2$ of a smooth algebraic curve $\G$
with marked points $p_\a$ can be defined by any pair of harmonic function $x_j(p)$ on $\G_0=\G\setminus\{p_\a\}$.

Let $x(p)$ be a harmonic function, then locally there exists a unique up to an additive constant conjugate harmonic function $y(p)$. Hence, $x(p)$ uniquely defines the differential $d\zeta=dx+idy$, which by construction is {\it imaginary normalized} holomorphic differential on $\G_0$. One can specify asymptotic behavior of $x(p)$ near the marked point by
the requirement that $d\zeta$ is meromorphic on $\G$ and has a fixed singular part at the marked points.

Recall that explicitly choosing a singular part of a pole of order at $n_\a+1$ near $p_\a$ means on a small neighborhood of $p_\a$ choosing: (i) a coordinate $z_\a$ such that $z_\a(p_\a)=0$; (ii) polynomial $R_\a$ of the form $R_\a=\sum_{i=0}^{n_\a}r_{\a,\,i}z_\a^{-i-1}$, and identifying pairs $(z_\a,R_\a)$ and $(w_\a, R'_\a)$ if $R'\,dw_\a=R_\a\,dz_\a+O(1)\,dz_\a,\, \, w_\a=w_\a(z_\a)$. The coefficient $r_{\a,\,0}$ is the residue of the singular part, i.e. for singular parts with no residues $r_{\a,\,0}=0$.

Non-degeneracy of the imaginary part of Riemann matrix of $b$-periods of normalized holomorphic differential on a smooth genus $g$ algebraic curve $\G$ implies that:

\begin{lem}\label{1} For any fixed singular parts of poles with pure real residues, there exists a unique meromorphic differential $\Psi$, having prescribed singular part at $P_\a$ and such that all its periods on $\G$ are imaginary, i.e.
\beq\label{realnorm}
{\rm Re}\, \left(\oint_c d\zeta\right)=0, \ \ \forall\  c\in H^1(\G,\mathbb Z)
\eeq
\end{lem}
\noindent
(for detailed proof see Proposition 3.4 in \cite{kr-grush1}).

The generalized amoebas considered above correspond to the case when $d\zeta$ has simple poles at the marked points. As it seems for the author the following two examples show that properties of maps $\chi$ defined by imaginary normalized differentials with higher order poles deserve a systematic study.

\bigskip

{\it Example 1.} Let $d\zeta_1, d\zeta_2$ be imaginary normalized differentials on $\G$ having pole at a marked point $p_0$ of the form $(z^{-2}+O(1))dz$  and $i(z^{-2}+O(1))dz$. Notice, that a different choice of the local coordinate $z$ corresponds to a linear transformation of the pair of differentials.

At first glance a notion of the corresponding amoeba is trivial. Indeed, it is easy to see that the map $\chi$ is one-to-one in the neighborhood of $p_0$. Hence, $\chi$ has degree $1$, and therefore, the image of $\chi$ is $\bbR^2$. Non-trivial nature of $\chi$ is reflected by a compact set in $\bbR^2$, which is the complement to the set whose points have just one preimage on $\G_0$. It can be checked that for an elliptic curve this complement set is bounded by piecewise concave curve with four cusps. An attempt to visualize the map $\chi$ the elliptic curve $\G_0$ reveals 3D creature which looks like jellyfish
(continuing "biological" terminology by GKZ).

\bigskip

{\it Example 2.} The following example of a pair of imaginary normalized differentials having poles of the form $d\zeta_1=(-z^-2+O(1))$ and $d\zeta_2=(-2z^{-3}+O(1))dz)$ is connected with the spectral theory of nonstationary Shr\"odinger equation (see details in \cite{kr-real},\cite{kr-spec}). It is easy to see that the corresponding map $\chi$ is of degree zero. There is one infinite connected complement of the image of $\chi$ which is bounded by a curve which is asymptotically is the parabola $x_2=x_1^2$. As shown in \cite{kr-spec} for the case of $M$-curves and one puncture fixed under anti-involution $\tau$ the map $\chi$ is $2:1$ outside of images of fixed ovals, which are boundaries of compact connected componets of $\A^c$. The gradient map $\nabla \rho$ restricted to $\G^+$ is one-to-one with the upper half plane of $\bbR^2$ with $g$ points removed.

\begin{rem}
The computations in the proof of Lemma \ref{critical} are pure local and can be summarized as follows:

Let $x_1(p)$ and $x_2(p)$ be a pair of harmonic functions in the domain $p\in\D$. They consecutively define: a pair of conjugate harmonic functions $y_j(p)$; the imaginary normalized holomorphic differentials $d\zeta_j=dx_j+idy_j$, and the function $R(p)$ (given by (\ref{R})).

\begin{cor} In the subdomain $\D_0:= \{p\in \D_0|\, \im R(p)\neq 0\}$, where the harmonic functions $x_j(p)$ define a system of local coordinates, there is a unique convex function $G(p)$ whose gradient equals

$$\nabla G(x_1,x_2)={\rm sgn}\, ( \im R) \left(\begin{array}{c} -y_2(x_1,x_2)\\
\ \,y_1(x_1,x_2) \end{array} \right).$$

\end{cor}

The generating function $G$ is a building block of the generalized Ronkin function which
at noncritical points $x$ of the image of the harmonic map $\chi:\D\longmapsto \bbR^2$ is equal to
$$\rho(x):=\sum_{p\in \chi^{-1}(x)} G(p).$$

\end{rem}

\end{document}